\documentclass[titlepage,11pt]{article}
\usepackage{latexsym}
\usepackage{amsfonts}
\usepackage{amsmath}
\usepackage{mathtools}
\usepackage{tikz}
\newcommand\blackslug{\hbox{\hskip 1pt \vrule width 4pt height 8pt depth 1.5pt
        \hskip 1pt}}
\newcommand\bbox{\hfill \quad \blackslug \bigbreak}

\def\LL{,\ldots,}

\newcommand{\mac}{\mathcal}

\title{A counterexample to the coarse Menger conjecture}
\author{
Tung Nguyen\thanks{Supported by AFOSR grants
A9550-19-1-0187 and FA9550-22-1-0234, and by NSF grants  DMS-1800053 and DMS-2154169.}\\
Princeton University,\\ Princeton, NJ 08544, USA
\and
Alex Scott\thanks{Supported by EPSRC grant EP/X013642/1}\\
University of Oxford, \\
Oxford, UK
\and
Paul Seymour\thanks{Supported by AFOSR grants
A9550-19-1-0187 and FA9550-22-1-0234, and by NSF grants  DMS-1800053 and DMS-2154169.}\\
Princeton University,\\ Princeton, NJ 08544, USA}

\date{November 19, 2023; revised \today}

\newtheorem{thm}{}[section]

\newcommand{\Proof}{\noindent{\bf Proof.}\ \ }

\begin{document}
\maketitle
\begin{abstract}
Menger's well-known theorem from 1927 characterizes when it is possible to find $k$ vertex-disjoint paths between two sets of vertices in a graph $G$.  Recently, Georgakopoulos and Papasoglu and, independently,
Albrechtsen, Huynh, Jacobs, Knappe and Wollan conjectured a coarse analogue of Menger's theorem, when the $k$ paths are required to be pairwise at some distance at least $d$.
The result is known for $k\le 2$, but we will show that it is false for all $k\ge 3$, even if $G$ is constrained to have maximum degree
at most three. We also give a simpler proof of the result when $k=2$.

\end{abstract}

\section{Introduction}

One of the foundational results in graph theory is a theorem of
Menger \cite{menger} from 1927 on vertex-disjoint paths:

\begin{thm}\label{menger}
{\bf Menger's theorem:}
Let $k\ge1$ be an integer, let $G$ be a graph and let $S,T\subseteq V(G)$.  Then either 
\begin{itemize}
\item there are $k$ vertex-disjoint paths between $S$ and $T$; or 
\item there is a set $X$ of at most $k-1$ vertices such that
every path between $S$ and $T$ contains a vertex of $X$.
\end{itemize}
\end{thm}

Recently, there has been significant interest in versions of Menger's theorem where the paths are not too close together: in other words, when are looking for $k$ paths between $S$, $T$ that are pairwise at distance at least $d$.\footnote{If $X$ and $Y$ are sets of vertices (or subgraphs) of a graph $G$, then $d(X,Y)$ denotes the distance between $X,Y$, that is, the number of edges in the shortest path of $G$ with one end in $X$ and the other in $Y$.}
One motivation for this comes from the developing area of `coarse graph theory', which is concerned with the large-scale geometric structure of graphs.  A central goal here is to find coarse analogues of classical graph-theoretic results (see Georgakopoulos and Papasoglu~\cite{agelos}).  A second motivation comes from work on the complexity of finding paths and other structures at distance at least 
two, with a substantial line of research starting with Bienstock \cite{bienstock} (see also Kawarabayashi and Kobayashi \cite{kk}; and a recent paper of Balig\'acs and J. MacManus \cite{bm} for larger distances).

A coarse analogue of Menger's theorem was conjectured by
Albrechtsen, Huynh, Jacobs, Knappe and Wollan~\cite{wollan}, and independently by Georgakopoulos and Papasoglu~\cite{agelos}:
\begin{thm}\label{conj}
{\bf Coarse Menger Conjecture:} For all integers $k,d\ge 1$ there exists $\ell>0$ with the following property.
Let $G$ be a graph and let $S,T\subseteq V(G)$.  Then either 
\begin{itemize}
\item there are $k$ paths between $S,T$, pairwise at distance at least $d$; or 
\item there is a set $X\subseteq V(G)$
with $|X|\le k-1$ such that every path between $S,T$ contains a vertex with distance at most $\ell$ from some member of $X$.
\end{itemize}
\end{thm}
Both sets of authors proved the result for $k=2$, but the $k\ge 3$ case remained open. 
The case $d=3$ is of special interest, because it is easy to see that if the result is true when $d=3$ (for some value of $k$) then it is true for all 
$d\ge 3$ and the same value of $k$ (apply the result when $d=3$ to the $d$th power of $G$).

Here we give a counterexample with
$d=k=3$, and with maximum degree at most three (taking multiple copies, or disjoint unions with paths, gives counterexamples for all $k,d\ge3$): 

\begin{thm}\label{construction}
For every $\ell>0$, there exists a graph $G$ with maximum degree 3 and vertex sets $S,T\subseteq V(G)$ such that 
\begin{itemize}
\item $G$ does not contain three paths between $S,T$, pairwise at distance at least three; and 
\item there is no set $X\subseteq V(G)$
with $|X|\le 2$ such that every path between $S,T$ contains a vertex with distance at most $\ell$ from some member of $X$.
\end{itemize}
\end{thm}

The rest of the paper is structured as follows.  In the next section, we give our construction proving \ref{construction}.  In section \ref{simpler}, we give a simpler proof of the $k=2$ case.  We conclude in section \ref{conclusion} with further discussion about possible variants of \ref{conj}.


\section{The counterexample}

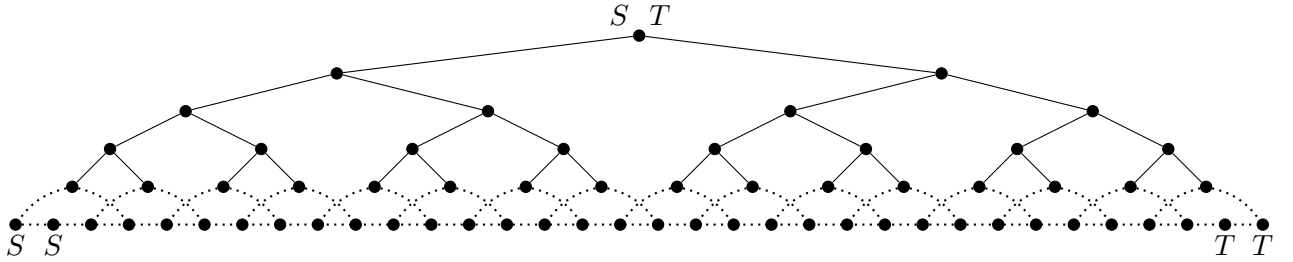
\begin{figure}[h!]
\centering

\begin{tikzpicture}[scale=1/2,auto=left]

\tikzstyle{every node}=[inner sep=1.5pt, fill=black,circle,draw]
\node (v1) at (1,0) {};
\node (v2) at (2,0) {};
\node (v3) at (3,0) {};
\node (v4) at (4,0) {};
\node (v5) at (5,0) {};
\node (v6) at (6,0) {};
\node (v7) at (7,0) {};
\node (v8) at (8,0) {};
\node (v9) at (9,0) {};
\node (v10) at (10,0) {};
\node (v11) at (11,0) {};
\node (v12) at (12,0) {};
\node (v13) at (13,0) {};
\node (v14) at (14,0) {};
\node (v15) at (15,0) {};
\node (v16) at (16,0) {};
\node (v17) at (17,0) {};
\node (v18) at (18,0) {};
\node (v19) at (19,0) {};
\node (v20) at (20,0) {};
\node (v21) at (21,0) {};
\node (v22) at (22,0) {};
\node (v23) at (23,0) {};
\node (v24) at (24,0) {};
\node (v25) at (25,0) {};
\node (v26) at (26,0) {};
\node (v27) at (27,0) {};
\node (v28) at (28,0) {};
\node (v29) at (29,0) {};
\node (v30) at (30,0) {};
\node (v31) at (31,0) {};
\node (v32) at (32,0) {};
\node (v33) at (33,0) {};
\node (v34) at (34,0) {};

\draw[dotted, thick] (v1)--(v34);

\node (u2) at (2.5,1) {};
\node (u4) at (4.5,1) {};
\node (u6) at (6.5,1) {};
\node (u8) at (8.5,1) {};
\node (u10) at (10.5,1) {};
\node (u12) at (12.5,1) {};
\node (u14) at (14.5,1) {};
\node (u16) at (16.5,1) {};
\node (u18) at (18.5,1) {};
\node (u20) at (20.5,1) {};
\node (u22) at (22.5,1) {};
\node (u24) at (24.5,1) {};
\node (u26) at (26.5,1) {};
\node (u28) at (28.5,1) {};
\node (u30) at (30.5,1) {};
\node (u32) at (32.5,1) {};

\draw[dotted,thick] (u2) to [bend right=20] (v1);
\draw[dotted,thick] (u4) to [bend right=20] (v3);
\draw[dotted,thick] (u6) to [bend right=20] (v5);
\draw[dotted,thick] (u8) to [bend right=20] (v7);
\draw[dotted,thick] (u10) to [bend right=20] (v9);
\draw[dotted,thick] (u12) to [bend right=20] (v11);
\draw[dotted,thick] (u14) to [bend right=20] (v13);
\draw[dotted,thick] (u16) to [bend right=20] (v15);
\draw[dotted,thick] (u18) to [bend right=20] (v17);
\draw[dotted,thick] (u20) to [bend right=20] (v19);
\draw[dotted,thick] (u22) to [bend right=20] (v21);
\draw[dotted,thick] (u24) to [bend right=20] (v23);
\draw[dotted,thick] (u26) to [bend right=20] (v25);
\draw[dotted,thick] (u28) to [bend right=20] (v27);
\draw[dotted,thick] (u30) to [bend right=20] (v29);
\draw[dotted,thick] (u32) to [bend right=20] (v31);

\draw[dotted,thick] (u2) to [bend left=20] (v4);
\draw[dotted,thick] (u4) to [bend left=20] (v6);
\draw[dotted,thick] (u6) to [bend left=20] (v8);
\draw[dotted,thick] (u8) to [bend left=20] (v10);
\draw[dotted,thick] (u10) to [bend left=20] (v12);
\draw[dotted,thick] (u12) to [bend left=20] (v14);
\draw[dotted,thick] (u14) to [bend left=20] (v16);
\draw[dotted,thick] (u16) to [bend left=20] (v18);
\draw[dotted,thick] (u18) to [bend left=20] (v20);
\draw[dotted,thick] (u20) to [bend left=20] (v22);
\draw[dotted,thick] (u22) to [bend left=20] (v24);
\draw[dotted,thick] (u24) to [bend left=20] (v26);
\draw[dotted,thick] (u26) to [bend left=20] (v28);
\draw[dotted,thick] (u28) to [bend left=20] (v30);
\draw[dotted,thick] (u30) to [bend left=20] (v32);
\draw[dotted,thick] (u32) to [bend left=20] (v34);

\node (t3) at (3.5,2) {};
\node (t7) at (7.5,2) {};
\node (t11) at (11.5,2) {};
\node (t15) at (15.5,2) {};
\node (t19) at (19.5,2) {};
\node (t23) at (23.5,2) {};
\node (t27) at (27.5,2) {};
\node (t31) at (31.5,2) {};

\draw (u2) -- (t3)--(u4);
\draw (u6) -- (t7)--(u8);
\draw (u10) -- (t11)--(u12);
\draw (u14) -- (t15)--(u16);
\draw (u18) -- (t19)--(u20);
\draw (u22) -- (t23)--(u24);
\draw (u26) -- (t27)--(u28);
\draw (u30) -- (t31)--(u32);

\node (s5) at (5.5,3) {};
\node (s13) at (13.5,3) {};
\node (s21) at (21.5,3) {};
\node (s29) at (29.5,3) {};

\draw (t3) -- (s5)--(t7);
\draw (t11) -- (s13)--(t15);
\draw (t19) -- (s21)--(t23);
\draw (t27) -- (s29)--(t31);

\node (r9) at (9.5,4) {};
\node (r25) at (25.5,4) {};

\draw (s5) -- (r9)--(s13);
\draw (s21) -- (r25)--(s29);

\node (q17) at (17.5,5) {};
\draw (r9) -- (q17)--(r25);

\tikzstyle{every node}=[]
\draw[above left] (q17) node []           {$S$};
\draw[above right] (q17) node []           {$T$};
\draw[below] (v1) node []           {$S$};
\draw[below] (v2) node []           {$S$};
\draw[below] (v33) node []           {$T$};
\draw[below] (v34) node []           {$T$};

\end{tikzpicture}

\caption{The dotted curves represent long paths.} \label{fig:counterexample}
\end{figure}

In this section, for each value of $\ell>0$,  we give an example of a graph $G$, and two subsets $S,T$ of $V(G)$,
such that there do not exist three paths between $S,T$, pairwise at distance at least three, and for every $X\subseteq V(G)$
with $|X|\le 2$, there is a path $P$ between $S,T$ such that its distance from $X$ is more than $\ell$. It is illustrated 
in figure~\ref{fig:counterexample}.  The two sets $S,T$ both have size three, and there is a vertex in $S\cap T$.

To show that this is a counterexample, we need to check that for every $X\subseteq V(G)$
with $|X|\le 2$, there is a path $P$ between $S,T$ such that its distance from $X$ is more than $\ell$; and that
there do not exist three paths between $S,T$, pairwise at distance at least three. The first is easy, so let us do it now.

\begin{thm}\label{countereasy}
Let $\ell\ge 1$ be an integer, and let $G$ be as in figure \ref{fig:counterexample}, where the binary tree has depth more than $2\ell+2$, and each of the dotted curves represents a path of length more than $2\ell$. Let $S,T$ be as shown in the figure.
If $X\subseteq V(G)$ with $|X|\le 2$, then there is a path $P$ between $S,T$ such that $d(X,P)>\ell$.
\end{thm}
\Proof Let $r$ be the root of the binary tree (the vertex at the top of the figure).
Since $r\in S\cap T$, we may assume there exists $x_1\in X$ such that its distance from 
$r$ is at most $\ell$; and since
all of the paths represented by the dotted curves (let us call them ``dotted paths'') have distance at least $2\ell+1$ from $r$ (because the binary tree has depth more than
$2\ell+2$), they all have distance more than $\ell$ from $x_1$.
Since
the path ($M$ say) at the bottom of the figure is between $S,T$, we may assume that
there exists $x_2\in X$ with distance at most $\ell$ from $M$. The set of vertices of $G$ with distance at most $\ell$ from $x_2$ 
is either a subset of the vertex set of one of the dotted paths in $M$, or it contains exactly one end $v$ of a dotted path of $M$,
and consists of $v$ together with  
subsets of the interiors of the (at most three) dotted paths incident with $v$. In either case, there is an end $v$ of one of the dotted paths in $M$, such that the set of vertices of $G$ with distance at most $\ell$ from $x_2$ is a subset of the set consisting of $v$
together with the interiors of the dotted paths incident with $v$. But for every choice of $v$, there is a path between $S,T$ 
made by a union of dotted paths, none of them incident with $v$; and consequently this path has distance more than $\ell$ from $X$.
This proves \ref{countereasy}.~\bbox

To prove the second statement, that there do not exist three paths between $S,T$, pairwise at distance at least $d$, we will prove
something stronger, by induction on the depth of the binary tree; and for that we need to allow the binary tree to have
small depth, and we need to set up some notation, so let us define the graph more carefully.

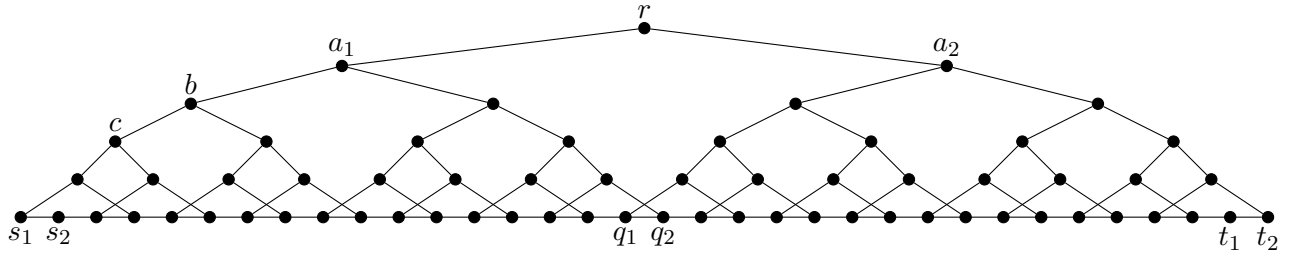
\begin{figure}[h!]
\centering

\begin{tikzpicture}[scale=1/2,auto=left]

\tikzstyle{every node}=[inner sep=1.5pt, fill=black,circle,draw]
\node (v1) at (1,0) {};
\node (v2) at (2,0) {};
\node (v3) at (3,0) {};
\node (v4) at (4,0) {};
\node (v5) at (5,0) {};
\node (v6) at (6,0) {};
\node (v7) at (7,0) {};
\node (v8) at (8,0) {};
\node (v9) at (9,0) {};
\node (v10) at (10,0) {};
\node (v11) at (11,0) {};
\node (v12) at (12,0) {};
\node (v13) at (13,0) {};
\node (v14) at (14,0) {};
\node (v15) at (15,0) {};
\node (v16) at (16,0) {};
\node (v17) at (17,0) {};
\node (v18) at (18,0) {};
\node (v19) at (19,0) {};
\node (v20) at (20,0) {};
\node (v21) at (21,0) {};
\node (v22) at (22,0) {};
\node (v23) at (23,0) {};
\node (v24) at (24,0) {};
\node (v25) at (25,0) {};
\node (v26) at (26,0) {};
\node (v27) at (27,0) {};
\node (v28) at (28,0) {};
\node (v29) at (29,0) {};
\node (v30) at (30,0) {};
\node (v31) at (31,0) {};
\node (v32) at (32,0) {};
\node (v33) at (33,0) {};
\node (v34) at (34,0) {};

\draw[] (v1)--(v34);

\node (u2) at (2.5,1) {};
\node (u4) at (4.5,1) {};
\node (u6) at (6.5,1) {};
\node (u8) at (8.5,1) {};
\node (u10) at (10.5,1) {};
\node (u12) at (12.5,1) {};
\node (u14) at (14.5,1) {};
\node (u16) at (16.5,1) {};
\node (u18) at (18.5,1) {};
\node (u20) at (20.5,1) {};
\node (u22) at (22.5,1) {};
\node (u24) at (24.5,1) {};
\node (u26) at (26.5,1) {};
\node (u28) at (28.5,1) {};
\node (u30) at (30.5,1) {};
\node (u32) at (32.5,1) {};

\draw[] (u2) to (v1);
\draw[] (u4) to  (v3);
\draw[] (u6) to  (v5);
\draw[] (u8) to  (v7);
\draw[] (u10) to (v9);
\draw[] (u12) to (v11);
\draw[] (u14) to (v13);
\draw[] (u16) to (v15);
\draw[] (u18) to (v17);
\draw[] (u20) to (v19);
\draw[] (u22) to (v21);
\draw[] (u24) to (v23);
\draw[] (u26) to (v25);
\draw[] (u28) to (v27);
\draw[] (u30) to (v29);
\draw[] (u32) to (v31);

\draw[] (u2) to (v4);
\draw[] (u4) to (v6);
\draw[] (u6) to (v8);
\draw[] (u8) to  (v10);
\draw[] (u10) to (v12);
\draw[] (u12) to (v14);
\draw[] (u14) to (v16);
\draw[] (u16) to (v18);
\draw[] (u18) to (v20);
\draw[] (u20) to (v22);
\draw[] (u22) to (v24);
\draw[] (u24) to (v26);
\draw[] (u26) to (v28);
\draw[] (u28) to (v30);
\draw[] (u30) to (v32);
\draw[] (u32) to (v34);

\node (t3) at (3.5,2) {};
\node (t7) at (7.5,2) {};
\node (t11) at (11.5,2) {};
\node (t15) at (15.5,2) {};
\node (t19) at (19.5,2) {};
\node (t23) at (23.5,2) {};
\node (t27) at (27.5,2) {};
\node (t31) at (31.5,2) {};

\draw (u2) -- (t3)--(u4);
\draw (u6) -- (t7)--(u8);
\draw (u10) -- (t11)--(u12);
\draw (u14) -- (t15)--(u16);
\draw (u18) -- (t19)--(u20);
\draw (u22) -- (t23)--(u24);
\draw (u26) -- (t27)--(u28);
\draw (u30) -- (t31)--(u32);

\node (s5) at (5.5,3) {};
\node (s13) at (13.5,3) {};
\node (s21) at (21.5,3) {};
\node (s29) at (29.5,3) {};

\draw (t3) -- (s5)--(t7);
\draw (t11) -- (s13)--(t15);
\draw (t19) -- (s21)--(t23);
\draw (t27) -- (s29)--(t31);

\node (r9) at (9.5,4) {};
\node (r25) at (25.5,4) {};

\draw (s5) -- (r9)--(s13);
\draw (s21) -- (r25)--(s29);

\node (q17) at (17.5,5) {};
\draw (r9) -- (q17)--(r25);

\tikzstyle{every node}=[]
\draw[above] (q17) node []           {$r$};
\draw[below] (v1) node []           {$s_1$};
\draw[below] (v2) node []           {$s_2$};
\draw[below] (v33) node []           {$t_1$};
\draw[below] (v34) node []           {$t_2$};
\draw[below] (v17)  node []           {$q_1$};
\draw[below] (v18)  node []           {$q_2$};
\draw[above ] (r9) node []           {$a_1$};
\draw[above ] (r25) node []           {$a_2$};
\draw[above ] (s5) node []           {$b$};
\draw[above ] (t3) node []           {$c$};

\end{tikzpicture}

\caption{The graph $G_k$. (In this figure, $k=6$.) The vertex $b$ is some vertex on the path between $a_1$ and $s_1$, different from and nonadjacent to $s_1$; it need not be adjacent to $a_1$, and
it might equal $a_1$. Its child on that path is $c$.}\label{fig:counter2}
\end{figure}

Let $k\ge 2$ be an integer. Take a uniform binary tree $B$ with depth $k$. Thus, $B$ has a root, and has
$2^k-1$ vertices, and every path from the root to one of the leaves has exactly $k-1$ edges.
Now add two more vertices, and let $Z$ consist of the set of leaves of $B$ together with the two new vertices;
and add a path $M$ with vertex set $Z$,
as shown in figure \ref{fig:counter2}. We call the resulting graph $G_k$. (To get the graph of the counterexample, we need to replace the edges 
of $G_{2\ell+3}$ incident with vertices in $Z$ by long paths, but let us not do that yet.) For each vertex $v$ of $B$ that is not a leaf, 
there is a copy of some $G_h$ with root $v$, formed by $v$ and its descendants in $B$, and two extra vertices; we will
apply the inductive hypothesis to these smaller graphs. Let $M$ have ends $s_1,t_2$, and let $s_2,t_1$
be the neighbours in $M$ of $s_1,t_2$ respectively. 

We will show:
\begin{thm}\label{counterbetter}
For each integer $k\ge 2$, if $P,Q$ are vertex-disjoint paths of $G_k$ between $\{s_1,s_2\}$ and $\{t_1,t_2\}$, then either
\begin{itemize}
\item there is a path of length at most two between $V(P)$ and $V(Q)$ such that none of its edges have ends in $Z$, or
\item one of $P,Q$ is the path of $B$ between $s_1,t_2$.
\end{itemize}
\end{thm}
\Proof
We proceed by induction on $k$. If $k=1$ the result is clear, so we assume that $k\ge 2$ and the result holds for all smaller values
of $k$. Suppose that $P,Q$ are vertex-disjoint paths of $G_k$ between $\{s_1,s_2\}$ and $\{t_1,t_2\}$, and the first bullet of the 
theorem is false. Let $r$ be the root of $B$, let $a_1$ be the neighbour of $r$ in the path of $B$ between $r, s_1$, and define $a_2$ similarly.
Suppose first that neither of $P,Q$ contains $r$. Then we may assume that $P$ contains $q_1$ and $Q$ contains $q_2$, where $q_1,q_2$
are as shown in the figure; and so from the inductive hypothesis applied to the two copies of $G_{k-1}$ with roots $a_1$ and $a_2$,
$P$ contains the path of $B$ between $s_1,q_2$, and $Q$ contains the path of $B$ between $q_1,t_2$. So $a_1\in V(P)$ and $a_2\in V(Q)$;
but $r$ is adjacent to both $a_1,a_2$, a contradiction.

So $r$ belongs to one of $P,Q$, say to $P$. We need to show that $P$ contains the two subpaths of $B$ between $r,s_1$ and between $r,t_2$.
Suppose it does not contain the first, say; and let us call this path $A$.  Choose a maximal common subpath of $P,A$, with ends
$a_1, b$ say. (Possibly $b=a_1$.) We see that $b\ne s_1$, since $A \not\subseteq P$; and $b$ is not adjacent to $s_1$, since
one of $P,Q$ must use the edge of $B$ incident with $s_1$. Thus the subpath of $A$ between $b,s_1$ has length $h$ say, where $h\ge 2$.
Let $c$ be the child of $b$ in $V(A)$. It follows that neither of $P,Q$ contains $c$; and this contradicts
the inductive hypothesis, applied to the copy of $G_{h-1}$ with root $c$. This proves \ref{counterbetter}.~\bbox

From \ref{counterbetter}, it follows easily that the graph of figure \ref{fig:counterexample}, with $S,T$ as shown, has the property that there do not exist three paths between $S,T$, pairwise at distance at least three. Because if $P,Q,R$ are three such paths, one of them
contains the root $r$, and so consists just of the vertex $r$; and so the other two cannot use $r$, and so contradict \ref{counterbetter}.

\section{A simpler proof}\label{simpler}
As noted in the introduction, the case $k=2$ of the coarse Menger conjecture is true.
Georgakopoulos and Papasoglu~\cite{agelos} gave a proof in geometric language, while
Albrechtsen, Huynh, Jacobs, Knappe and Wollan~\cite{wollan} gave a relatively complicated combinatorial proof.  Here, we give a slightly simpler combinatorial proof.

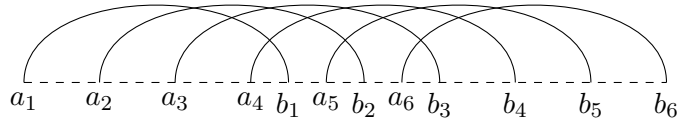
\begin{figure}[h!]
\centering

\begin{tikzpicture}[scale=1/2,auto=left]

\draw[dashed]  (0,0)--(17,0);

\draw (0,0) to [bend left = 90] (7,0);
\draw (2,0) to [bend left = 90] (9,0);
\draw (4,0) to [bend left = 90] (11,0);
\draw (6,0) to [bend left = 90] (13,0);
\draw (8,0) to [bend left = 90] (15,0);
\draw (10,0) to [bend left = 90] (17,0);

\tikzstyle{every node}=[]
\draw[below] (0,0) node []           {$a_1$};
\draw[below] (2,0) node []           {$a_2$};
\draw[below] (4,0) node []           {$a_3$};
\draw[below] (6,0) node []           {$a_4$};
\draw[below] (8,0) node []           {$a_5$};
\draw[below] (10,0) node []           {$a_6$};
\draw[below] (7,0) node []           {$b_1$};
\draw[below] (9,0) node []           {$b_2$};
\draw[below] (11,0) node []           {$b_3$};
\draw[below] (13,0) node []           {$b_4$};
\draw[below] (15,0) node []           {$b_5$};
\draw[below] (17,0) node []           {$b_6$};

\end{tikzpicture}

\caption{Intervals in standard form (the curves join the ends of the intervals).}\label{fig:standardform}
\end{figure}

We need some lemmas about intervals. An {\em interval} is a pair $(a,b)$ of integers with $a\le b$, and its {\em length} is $b-a$.
If $\mac H$ is a set of intervals that can be written $\mac H=\{(a_i,b_i):1\le i\le t\}$ for some $t\ge 1$, such that
$0\le a_1<a_2<\cdots <a_t\le n$, and $0\le  b_1<\cdots<b_t\le n$, we call this the {\em standard form} for $\mac H$. Thus $\mac H$ has a 
standard form if and only there do not exist distinct $(a,b),(c,d)\in \mac H$ with $a\le c\le d\le b$. 
(See figure \ref{fig:standardform}.)

If $(a,b), (c,d)$ are intervals, we say $(a,b)$ {\em captures} $(c,d)$ if $a\le c\le d\le b$; and a set $\mac H$ of intervals
{\em captures}
$(c,d)$ if there exists $(a,b)\in \mac H$ that captures $(c,d)$.
For integers $\ell, n$ with $0<\ell\le n$, we say a set $\mac H$ of intervals
is {\em $\ell$-powerful} in $(0,n)$ if
\begin{itemize}
\item $0\le a\le b\le n$ for all $(a,b)\in \mac H$; and
\item $\mac H$ captures every interval of length $\ell$ captured by $(0,n)$.
\end{itemize}

\begin{figure}[h!]
\centering

\begin{tikzpicture}[scale=1/2,auto=left]

\draw[dashed]  (0,0)--(17,0);

\draw (2,0) to [bend left = 90] (9,0);
\draw (4,0) to [bend left = 90] (11,0);
\draw (6,0) to [bend left = 90] (13,0);

\tikzstyle{every node}=[]
\draw[below] (2,0) node []           {$a_i$};
\draw[below] (4,0) node []           {$a_{i+1}$};
\draw[below] (6,0) node []           {$a_j$};
\draw[below] (9,0) node []           {$b_i$};
\draw[below] (11,0) node []           {$b_{i+1}$};
\draw[below] (13,0) node []           {$b_j$};

\end{tikzpicture}

\caption{For the proof of \ref{getdisjt}.}\label{fig:getdisjt}
\end{figure}
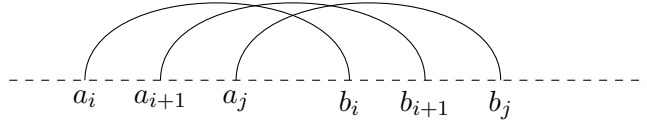

Our objective is, starting with an $\ell$-powerful set $\mac H$ of intervals, to find a subset $\mac H'$ of $\mac H$ which is still
$\ell'$-powerful
for some $\ell'$ not much smaller than $\ell$, such that all the ends of all the intervals in $\mac H'$ are far apart. (Not exactly:
in standard form, we cannot arrange that $a_{i+2}-b_i$ is large for each $i$, and it might even be negative,
and we cannot arrange that $a_2-a_1$ and $b_t-b_{t-1}$
are large, but we will make all the other differences large). We begin with:

\begin{thm}\label{getdisjt}
If $0<\ell\le n$ and $\mac H$ is a set of intervals that is $\ell$-powerful in $(0,n)$, and minimal with this property, then 
in standard form $a_j\ge b_{i}-\ell+2$ for all $i,j\in \{1\LL t\}$ with $j\ge i+2$.
\end{thm}
\Proof (See figure \ref{fig:getdisjt}.)

Let $1\le i,j\le t$ with $j\ge i+2$. From the minimality, there exists $h\in \{0,n-\ell\}$ such that $(a_{i+1},b_{i+1})$ captures $(h,h+\ell)$,
and no other member of $\mac H$ captures $(h,h+\ell)$. Since $(a_i,b_i)$ does not capture $(h,h+\ell)$, and $a_i<a_{i+1}\le h$,
it follows that $b_i<h+\ell$; and since $(a_j,b_j)$ does not capture $(h,h+\ell)$, and $b_j\ge b_{i+1}\ge h+\ell$, it follows
that $a_j>h$. Hence $a_j\ge b_{i}-\ell+2$.
This proves \ref{getdisjt}.~\bbox

If $\mac H$ can be written in standard form, and $\mac H'\subseteq \mac H$, then $\mac H'$ can also be written in standard form. Some inequalities
about the corresponding numbers $a_i,b_i$ are {\em hereditary}; that is, if they hold for $\mac H$ then they also hold for all its subsets
$\mac H'$. 
For instance, the inequalities of \ref{getdisjt} are hereditary. 
Let $\mac H'\subseteq \mac H$, and $\mac H= \{(a_i,b_i):1\le i\le t\}$ and $\mac H' = \{(a_i',b_i'):1\le i\le t'\}$, in standard form.
If we know that $a_j\ge b_{i}-\ell+2$ for all $i,j\in \{1\LL t\}$ with $j\ge i+2$, then it follows that 
$a_j'\ge b_{i}'-\ell+2$ for all $i,j\in \{1\LL t'\}$ with $j\ge i+2$. Most of the other inequalities we will prove are hereditary,
and we leave verifying this to the reader.

\begin{thm}\label{int2}
If $0<2\ell\le n$ and $\mac H$ is a set of intervals that is $2\ell$-powerful in $(0,n)$, 
then there exists $\mac H'\subseteq \mac H$, minimally $\ell$-powerful in $(0,n)$, 
that can be written $\mac H'=\{(a_i,b_i):1\le i\le t\}$ in standard form, such that
\begin{itemize}
\item $a_j\ge b_{i}-\ell+2$ for $1\le i,j\le t$ with $j\ge i+2$; and
\item $b_i-b_{i-1}\ge \ell $ for $1< i< t$, and $b_{i-1}-a_i\ge \ell$ for $1< i\le  t$.
\end{itemize}
\end{thm}
\Proof
We may assume that $\mac H$  is minimally $2\ell$-powerful in $(0,n)$.
We define a sequence $(a_i,b_i) (1\le i\le t)$ of members of $\mac H$ inductively as follows.
Let $(a_1,b_1)$ be the (unique) member of $\mac H$ that captures $(0,2\ell)$. Inductively, having defined
$(a_i,b_i)$:
\begin{itemize}
\item If $b_i < n - \ell$, choose $(a_{i+1}, b_{i+1}) \in \mac H$ 
capturing $(b_{i} - \ell, b_{i} + \ell)$ with $b_{i+1}$ as large as possible.
\item If $n - \ell \le b_{i} < n$, let $(a_{i+1}, b_{i+1})\in \mac H$ be the (unique) member of $\mac H$
capturing $(n - 2\ell, n)$.
\item If $b_i = n$, let $t=i$; the inductive definition is complete.
\end{itemize}
Let $\mac H_1 = \{(a_i,b_i):1\le i\le t\}$. From the construction, it follows that this expression of $H_1$ is in standard form.
Moreover, $b_i-b_{i-1}\ge \ell $ for $1< i< t$, and $b_{i-1}-a_i\ge \ell$ for $1< i\le  t$, that is, $H_1$ satisfies the second bullet of the theorem. 

We claim that $\mac H_1$ is $\ell$-powerful in $(0,n)$. To see this, let $0\le h\le n-\ell$; we must show that some member of $\mac H_1$
captures $(h,h+\ell)$. We may assume that $h<n-2\ell$, since $(a_t,b_t)$ captures $(n-2\ell, n)$.
Choose $i\in \{1\LL t\}$ maximal such that $a_i\le h$. We may assume that $b_i<h+\ell$, since otherwise
$(a_i,b_i)$ captures $(h,h+\ell)$. Hence $b_i<n-\ell$, and so $(a_{i+1},b_{i+1})$ captures $(b_{i}-\ell,b_{i}+\ell)$, from the construction.
From the choice of $i$, $a_{i+1}>h$; but $a_{i+1}\le b_i-\ell$, so $b_i\ge a_{i+1}+\ell>h+\ell$, a contradiction.
This proves that $\mac H_1$ is $\ell$-powerful in $(0,n)$. Choose $\mac H'\subseteq \mac H_2$, minimal such that it is $\ell$-powerful in $(0,n)$.
Then by \ref{getdisjt}, the first bullet of the theorem holds; and the second still holds, since it is hereditary.
This proves \ref{int2}.~\bbox

\begin{figure}[h!]
\centering

\begin{tikzpicture}[scale=1,auto=left]

\draw (0,0) to [bend left = 90] (4,0);
\draw (2,0) to [bend left = 90] (6,0);
\draw (4,0) to [bend left = 90] (8,0);
\draw (6,0) to [bend left = 90] (10,0);
\draw (8,0) to [bend left = 90] (12,0);
\draw (10,0) to [bend left = 90] (14,0);

\tikzstyle{every node}=[]
\draw[below] (0,0) node []           {$a_1$};
\draw[below] (2,0) node []           {$a_2$};
\draw[below] (4,0) node []           {$a_3\&b_1$};
\draw[below] (6,0) node []           {$a_4\&b_2$};
\draw[below] (10,0) node []           {$a_t\&b_{t-2}$};
\draw[below] (12,0) node []           {$b_{t-1}$};
\draw[below] (14,0) node []           {$b_t$};

\draw[white, fill=white] (6.8,0) rectangle (9.2,1.2);
\draw[dashed]  (0,0)--(14,0);

\draw[dotted,thick] (7,.6) -- (9,.6);
\end{tikzpicture}

\caption{The output of \ref{mainint}. $a_3, b_1$ have been drawn in the same place because we do not know which is larger.
}\label{fig:mainint}
\end{figure}
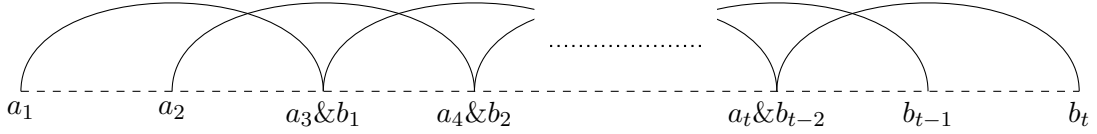

\begin{thm}\label{mainint}
If $0<4\ell\le n$ and $\mac H$ is a set of intervals that is $4\ell$-powerful in $(0,n)$,
then there exists $\mac H'\subseteq \mac H$, $\ell$-powerful in $(0,n)$,                        
which can be written $\mac H'=\{(a_i,b_i):1\le i\le t\}$ in standard form, such that
the order of the numbers $a_1\LL a_t$ and $b_1\LL b_t$
is:
\begin{align*}
&0=a_1< a_2<\min(a_3,b_1)<\max(a_3,b_1)<\min(a_4,b_2)<\max(a_4,b_2)<\cdots\\
&\cdots<\min(a_t,b_{t-2})<\max(a_t,b_{t-2})<b_{t-1}<b_t=n.
\end{align*}
(See figure \ref{fig:mainint}.)
Every two of $a_1\LL a_t,b_1\LL b_t$ differ by at least $\ell$, except possibly the pairs $(a_1,a_2)$, $(b_{t-1},b_t)$, and
$(a_i,b_{i-2})$ for $3\le i\le t$.
\end{thm}
\Proof
From \ref{int2} with $\ell$ replaced by $2\ell$, there exists $\mac H_1\subseteq \mac H$, $2\ell$-powerful in $(0,n)$, 
that satisfies the inequalities of \ref{int2}.
By \ref{int2} again, applied to $\mac H_1$, with the order of $\{0\LL n\}$ reversed, there exists $\mac H_2\subseteq \mac H_1$,
minimally $\ell$-powerful in $(0,n)$, such that in standard form,
\begin{itemize}
\item $b_i-b_{i-1}\ge 2\ell $ for $1< i< t$, and $b_{i-1}-a_i\ge 2\ell$ for $1< i\le  t$ (since these inequalities are true
for $\mac H_1$ and are hereditary);
\item $a_{i+1}-a_{i}\ge \ell$ for $1<i<t$; and
\item $a_j\ge b_i-\ell+2$ for all $i,j\in \{1\LL t\}$ with $j\ge i+2$.
\end{itemize}
It follows that $b_j-b_i\ge \ell$ for all $i,j\in \{1\LL t\}$ with $j>i$ and $(i,j)\ne (t-1,t)$;
that $a_j-a_i\ge \ell$ for all $i,j\in \{1\LL t\}$ with $j>i$ and $(i,j)\ne (1,2)$. Now let $i,j\in \{1\LL t\}$; we want to show that
$|a_j-b_i|\ge \ell$, unless $i=j-2$.
If $i>j$, then $b_i\ge b_j+\ell \ge a_j+\ell$ as required. If $j=i$, the minimality of $\mac H_2$ implies that $b_i-a_i\ge \ell$, as required.
If $i=j-1$, we have $b_{i}-a_j\ge 2\ell$ from the first bullet above, as required. Now suppose that $i\le j-3$. Hence
$b_{j-2}\le a_j+\ell-2$, and $b_{j-2}\ge b_{j-3}+2\ell$; so $a_j+\ell-2\ge b_{j-3}+2\ell$, that is,
$a_j\ge b_{j-3}+\ell+2\ge b_i+\ell$ as required. This proves \ref{mainint}.~\bbox

Now let us deduce \ref{conj} when $k=2$. As we said earlier,
it suffices to handle the case when $d=3$, so we will prove:

\begin{thm}\label{truecase}
Let $G$ be a non-null graph and let $S,T\subseteq V(G)$; then either
\begin{itemize}
\item there are two paths between $S,T$ with distance at least three; or
\item there exists $x\in V(G)$
such that every path between $S,T$ contains a vertex with distance at most $161$ from $x$.
\end{itemize}
\end{thm}
\Proof
Let $c\ge  7$ and $\ell\ge  2c+5$ be integers: we will prove that the theorem holds with $161$ replaced by $8\ell+c+2$.
We assume the second outcome is false.
We may assume that $G$ is connected, and
there is a path between $S,T$; let $R$ be such a path with minimum length, and let its vertices be $r_1\LL r_{n-1}$ in order,
where $r_1\in S$ and $r_{n-1}\in T$. There is a path $P$ between $S,T$ with $d(r_1,P)> 8\ell+c+2$, and we can assume that $d(P,R)\le 2$,
so $n-2\ge 8\ell+c+1$. In particular, $S\cap T=\emptyset$ and there are no edges between $S,T$.

Let $W$ be the set of all vertices with distance at most $c$ from $R$. We call the set of all vertices $v$ with $d(v,R)=c$ the {\em surface}. For each $v\in V(G)$,
if $v\in S$ and $d(v,R)\ge c+4$, let $a(v)=0$, and otherwise 
let $a(v)$ be the smallest $i\in \{1\LL n-1\}$ such that there is a path between $v$ and $r_i$ of length $d(v,R)$. Similarly, 
if $v\in T$ and $d(v,R)\ge c+4$ let $b(v)=n$, and otherwise let $b(v)$
be the largest $i\in \{1\LL n-1\}$ such that there is a path between $v$ and $r_i$ of length $d(v,R)$. 
For each $w\in W$, it follows from the choice of $R$ that $0\le b(w)-a(w)\le 2d(w,R)\le 2c$. 

Let $\mathcal{C}$ be the set of components of $G\setminus W$, and for each $C\in \mathcal{C}$, let $N(C)$ be the set of vertices
in $W$ with a neighbour in $V(C)$.  Thus $N(C)\ne \emptyset$, and $N(C)$ is a subset of the surface.
Let 
$a(C)$ be the smallest value of $a(v)$ for $v\in V(C)$, and let $b(C)$ be the largest value of $b(v)$ for $v\in V(C)$.
Thus $0\le a(C)\le b(C)\le n$ for each $C\in \mathcal{C}$.
\\
\\
(1) {\em The set of intervals $\{(a(C),b(C)):C\in \mathcal{C}\}$ is $16\ell$-powerful in $(0,n)$.}
\\
\\
Let $8\ell\le m\le n-8\ell$; we must show that there exists $C\in \mathcal{C}$ such that $(a(C),b(C))$ captures $(m-8\ell,m+8\ell)$. 
Let $B$
be the set of all vertices of $G$ with distance at most $8\ell+c+2$ from $r_{m}$. 
Let $W_1$ be the set of all $w\in W\setminus B$ with
$a(w)\le m$, and let $W_2$ be the set of all $w\in W\setminus B$ with $b(w)\ge m$. Thus $W_1\cup W_2=W\setminus B$. If $w\in W_1$, 
then since $d(w,r_m)>8\ell+c+2$
and $d(w,r_{a(w)})\le c$, it follows that $d(r_{a(w)},r_m)>8\ell+2$, and so $a(w)\le m-8\ell-3$; and similarly if $w\in W_2$
then $b(w)\ge m+8\ell+3$.
Thus, if $w_1\in W_1$ and $w_2\in W_2$, then 
$$b(w_2)-a(w_1)\ge (m+8\ell+3) - (m-8\ell-3) = 16\ell+6.$$
There is a path 
between $w_1$ and $S$ of length at most $a(w_1)+c-1$, and a path
between $w_2$ and $T$ of length at most $|E(R)|-b(w_2)+ c-1$,
and, since there is no path between $S,T$ of length less than $|E(R)|$, it follows that
\begin{align*}
d(w_1,w_2)&\ge |E(R)|- (a(w_1)+c-1)-(|E(R)|-b(w_2)+ c-1)\\
&=b(w_2)-a(w_1) -2(c-1)\ge (16\ell+6)-2(c-1)\ge 2.
\end{align*}
Consequently $W_1\setminus B$ and $W_2\setminus B$ are disjoint, and there are no edges between them. Moreover, if $w_1\in S$ and $w_2\in W_2$ then $b(w_2)\ge m+8\ell+3\ge 16\ell+3$, and similarly
$$d(w_1,w_2)\ge |E(R)|- (|E(R)|-b(w_2)+ c-1)\\
=b(w_2) -(c-1)\ge (16\ell+3)-2(c-1)\ge 2.$$
So $S\cap W_2=\emptyset$ and there are no edges between them, and the same for $W_1,T$.
Thus $W_1\cup S$ is disjoint from $W_2\cup T$ and there are no edges between them.

Since the second outcome of the theorem is false,
there is an $S-T$ path $P$ in $G$ disjoint from $B$.
Consequently, there is a subpath $P'$ of $P$
with first vertex in $W_1\cup S$, last vertex in $W_2\cup T$, and with no other vertex in $W\cup B$.
Since $W_1\cup S$ is disjoint from $W_2\cup T$ and there are no edges between them,
$P'$ has an internal vertex, which is therefore not in $W$.
But the only vertices of $P'$ in $W$ are ends of $P'$,
and so there exists
$C\in \mathcal{C}$ such that $V(P')\subseteq V(C)\cup W$.
We will show that $a(C)\le m-8\ell$ and $b(C)\ge m+8\ell$.
Let $w_1\in W_1\cup S$ and $w_2\in W_2\cup T$ be the ends of $P'$.
Either $w_1\in W_1$, or $w_1\in S\setminus W$ and $d(w_1,R)\ge c+4$, or $w_1\in S\setminus W$ and $d(w_1,R)\le c+3$, and we claim 
that $a(w_1) \le m-8\ell$ in each case.
In the first case, $w_1$ has a neighbour $v\in V(C)$; and since $d(v,R)=c+1=d(w_1,R)+1$, it follows that 
$a(C)\le a(v)\le a(w_1)\le m-8\ell-3\le  m-8\ell$ as required. In the other two cases, $w_1\in V(C)$, and so $a(C)\le a(w_1)$.
In the second case, $a(w_1)=0\le m-8\ell$ as required. In the third case, 
when $w_1\in S$ and $d(w_1,R)\le c+3$, it follows that $d(w_1,r_{a(w_1)})\le c+3$, and since $d(w_1, r_m)>8\ell+c+2$,
this implies that $a(w_1) \le m-8\ell$ as required. This proves that $a(C)\le m-8\ell$ in all three cases. 
Similarly $b(C)\ge m+8\ell$. This proves (1).

\bigskip

From \ref{mainint}, with $\ell$ replaced by $4\ell$, there exists $\mathcal{D}\subseteq \mathcal{C}$ such that the set of intervals 
$\{(a(C),b(C)):C\in \mathcal{D}\}$ is $4\ell$-powerful in $(0,n)$, and the members of $\mathcal{D}$ can be numbered $D_1\LL D_t$
such that:
\begin{itemize}
\item $0=a(D_1)<a(D_2)<\cdots <a(D_t)\le n$, and $0\le b(D_1)<b(D_2)<\cdots <b(D_t)=n$; and
\item $a(D_i)-b(D_{i-3})\ge 4\ell$ for $4\le i\le t$.
\end{itemize}
(We are not using the full strength of \ref{mainint} here; that will come later.)
\\
\\
(2) {\em For $1\le i,j\le t$, if $j\ge i+3$ then $d(D_i,D_j)\ge 4\ell-2c+2$.}
\\
\\
Suppose not. Hence there exist $w\in N(D_i)$ and $w'\in N(D_j)$ such that $d(w,w')<4\ell-2c$; and there is a path between 
$w, V(R)$
of length $c$, and its end in $R$ is some $r_k$ where $k\le b(D_i)$. Similarly 
there is a path between      
$w', V(R)$
of length $c$, and its end in $R$ is some $r_{k'}$ where $k'\ge a(D_j)$.
Thus there exist $k\le b(D_i)$ and $k'\ge a(D_j)$ such that there is a path between $r_k,r_{k'}$ of length less than $4\ell$.
Since $R$ is a shortest path, it follows that $4\ell-1\ge k'-k\ge a(D_j)-b(D_i)\ge 4\ell$, a contradiction. This proves (2).

\bigskip
Let $\Delta$ be the union of the vertex sets of the members of $\mathcal{D}$. If $X\subseteq V(G)$, let $\Delta(X)$
be the set of $D\in \mathcal{D}$ with $X\cap N(D)\ne \emptyset$. 
\\
\\
(4) {\em  If $X\subseteq V(G)$ is connected and $|\Delta(X)| \ge 4$, then $|X| \ge 30$.}
\\
\\
Since $|\Delta(X)|\ge 4$, we can choose
$i,j$ with $1\le i,j\le t$ and $j\ge i+3$ such that $D_i,D_j\in \Delta(X)$. Since $G[X]$
is connected, there is a path between $D_i,D_j$ of length at most $|X|+1$.
By (2), $d(D_i,D_j)\ge 4\ell-2c+2>30$, and so $|X|\ge 30$. This proves (4).

\bigskip

Let us say a {\em joint} is a subset $X\subseteq V(G)$, inducing a connected subgraph, such that 
every vertex of $X\cap W$ has distance at most $(|X|-1)/2$ from the surface, and 
either
\begin{itemize}
\item $|\Delta(X)|\ge2$ and $|X|\le 3$; or
\item $|\Delta(X)|\ge  3$ and $|X|\le 8$.
\end{itemize}
\noindent(4) {\em For every joint $X$, $d(X,R)\ge c-(|X|-1)/2$.}
\\
\\
If $v\in X$ and $d(v,R)\le c$, then $v\in W$ and so 
there is a path from $v$ to the surface of length at most $(|X|-1)/2$. Since every path from the surface to $R$ has length at least $c$,
it follows that $d(v,R)\ge c-(|X|-1)/2$.
This proves (4).

\bigskip
Let $Z$ be the union of all the joints. Each component
of $G[Z\cup \Delta]$ is called a {\em supercomponent}. Each supercomponent includes at least one member of $\mathcal{D}$, but may be
composed of many joints and members of $\mathcal{D}$.
\\
\\
(5) {\em Let $F_1,F_2$ be distinct supercomponents, and let $Q$ be a path of $G$ with ends $f_1\in V(F_1)$ and $f_2\in V(F_2)$.
\begin{itemize}
\item If both $f_1,f_2$ belong to joints then   $|Q|\ge 16
$; 
\item if exactly one of $f_1,f_2$ belongs to a joint $X$ then $|Q|\ge 8
$, 
and $|Q|\ge 24
$ if $|\Delta(X)|\ge 3$;
and 
\item if neither of $f_1,f_2$ belong to joints then $|Q|\ge 6$.
\end{itemize}
In any case, $d(F_1,F_2)\ge 5$.}
\\
\\
Let $Q^*$ be the set of vertices in the interior of $Q$. 
Suppose first that $f_i$ belongs to a joint $X_i$ for $i = 1,2$. Since
$|\Delta(X_1\cup X_2\cup Q^*)|\ge 4$ (since $|\Delta(X_i)|\ge 2$ for $i = 1,2$, 
and $\Delta(X_1)\cap \Delta(X_2)=\emptyset$), it follows from (4) that  $|X_1\cup X_2\cup Q^*|\ge 30$.
Since $|X_1|,|X_2|\le 8$, it follows that $|Q^*|\ge 14$, and so $|Q|\ge 16$ as claimed. Thus we may assume that $f_2\in D'$ for some $D'\in \mathcal{D}$ included in $F_2$,
and so $Q^*\cap N(D')\ne \emptyset$, and hence $D'\in \Delta(Q^*)$. 
Next, suppose that $f_1\in X_1$ for some joint $X_1$. Thus $f_1$ has distance at most $(|X_1|-1)/2$ from the surface, 
and so each vertex of $X_1\cup Q^*$ has distance at most 
$$(|Q|-2+(|X_1|-1)/2)/2\le (|X_1\cup Q^*|-1)/2$$ 
from the surface.
But $|\Delta(X_1\cup Q^*)|\ge \Delta(X_1)+1\ge 3$, and since $X_1\cup Q^*$ is not a joint (because $F_1,F_2$
are distinct supercomponents), it follows that $|X_1\cup Q^*|\ge 9$.
If  $|\Delta(X_1)|=2$, then $|X_1|\le 3$, and so $|Q^*|\ge 6$ and $|Q|\ge 8$; and if $|\Delta(X_1)|\ge 3$ then $|X_1|\le 8$ (because $X_1$
is a joint), and $|X_1\cup Q^*|\ge 30$
by (4), and so 
$|Q^*|\ge 22$ and $|Q|\ge 24$, as claimed. Finally, we may assume that 
$f_1\in D$ for some $D\in \mathcal{D}$ included in $F_1$, and so $D,D'\in \Delta(Q^*)$. Since every vertex of $Q^*\cap W$ has 
distance at most $(|Q^*|-1)/2$ from the surface, and $Q^*$ is not a joint, it follows that 
$|Q^*|\ge 4$ and so $|Q|\ge 6$. This proves (5).
\\
\\
(6) {\em If $F_1,F_2$ are distinct supercomponents, and $A$ is a path of length at most $c+1$
between $F_2$ and $R$, then 
$d(F_1,A)\ge 3$.}
\\
\\
Suppose that $f_1\in V(F_1)$ and $v\in V(A)$ have distance at most two. Let $f_2$ be the end of $A$ in $V(F_2)$.
Suppose first that
$f_1\in \Delta$. Consequently every path from $f_1$ to $V(R)$ has length more than $c$, and so the subpath of $A$ between $v$ and $R$
has length at least $c-1$; and therefore the subpath from $v$ to $f_2$ has length at most two.
Hence there is a path $Q$ of length at most four 
between $f_1,f_2$; so $|Q|\le 5$, contrary to (5).
This proves that $f_1\notin \Delta$, and so $f_1\in X_1$ for some joint $X_1$ of $F_1$. 
By (4), $d(f_1,R)\ge c-(|X_1|-1)/2$, and so
the subpath of $A$ between $v, R$ has length at least $c-(2+(|X_1|-1)/2)$. Hence the subpath of $A$ between $v,f_2$ has length 
at most $3+ (|X_1|-1)/2$. Consequently there is a path $Q$ between $f_1,f_2$ of length at most $5+ (|X_1|-1)/2$, and so 
$|Q|\le 6 + (|X_1|-1)/2$.
Since $|X_1|\le 8$, and
so $|Q|\le 9$, it follows from (5) that $f_2$ is not in a joint and $|\Delta(X_1)|=2$.
But then $|X_1|\le 3$, so $|Q|\le 7$, contrary to (5). This proves (6).


\bigskip

Let $\mathcal{F}$ be the set of all supercomponents. For each $F\in \mathcal{F}$, let $a(F)$ be the minimum of $a(v)$ over all $v\in V(F)$,
 and define $b(F)$ similarly.
Since for every $D\in \mathcal{D}$, there exists $F\in \mathcal{F}$
with $D\subseteq F$, and hence with $a(F)\le a(D)\le b(D)\le b(F)$, it follows that the set of intervals 
$\{(a(F),b(F)):F\in \mathcal{F}\}$ is $4\ell$-powerful. By \ref{mainint}, there exists $\mathcal{H}\subseteq \mathcal{F}$ such that 
$\{(a(F),b(F)):F\in \mathcal{H}\}$ is $\ell$-powerful, and it can be numbered as $\{H_1\LL H_s\}$,
where, writing $a(H_i) = a_i$ and $b(H_i)=b_i$ for each $i$:
\\
\\
(7) {\em The order of the numbers $a_1\LL a_s$ and $b_1\LL b_s$
is:
\begin{align*}
&0=a_1< a_2<\min(a_3,b_1)<\max(a_3,b_1)<\min(a_4,b_2)<\max(a_4,b_2)<\cdots\\
&\cdots<\min(a_s,b_{s-2})<\max(a_s,b_{s-2})<b_{s-1}<b_s=n.
\end{align*}
Every two of $a_1\LL a_s,b_1\LL b_s$ differ by at least $\ell$, except possibly the pairs $(a_1,a_2)$, $(b_{s-1},b_s)$, and
$(a_i,b_{i-2})$ for $3\le i\le s$.}

\bigskip

\begin{figure}[h!]
\centering

\begin{tikzpicture}[scale=1/2,auto=left]

\def\c{7}
\draw[dotted]  (0,\c)--(30,\c);
\tikzstyle{every node}=[inner sep=1.5pt, fill=black,circle,draw]
\node (r1) at (0,0) {};
\node (ra2) at (1.5,0) {};
\node (rb1) at (5,0) {};
\node (ra3) at (6.5,0) {};
\node (ra4) at (10,0) {};
\node (rb2) at (11.5,0) {};
\node (ra5) at (15,0) {};
\node (rb3) at (16.5,0) {};
\node (rb4) at (20,0) {};
\node (ra6) at (21.5,0) {};
\node (ra7) at (25,0) {};
\node (rb5) at (26.5,0) {};

\node (a1) at (0,\c+1) {};
\node (a2) at (2.5,\c) {};
\node (b1) at (5,\c) {};
\node (a3) at (7.5,\c) {};
\node (b2) at (10,\c) {};
\node (a4) at (12.5,\c) {};
\node (b3) at (15,\c) {};
\node (a5) at (17.5,\c) {};
\node (b4) at (20,\c) {};
\node (a6) at (22.5,\c) {};
\node (b5) at (25,\c) {};
\node (a7) at (27.5,\c) {};
\node (b6) at (30,\c) {};

\tikzstyle{every node}=[]
\draw[below] (r1) node []           {$r_1$};
\draw[below] (ra2) node []           {$r_{a_2}$};
\draw[below] (rb1) node []           {$r_{b_1}$};
\draw[below] (ra3) node []           {$r_{a_3}$};
\draw[below] (ra4) node []           {$r_{a_4}$};
\draw[below] (rb2) node []           {$r_{b_2}$};
\draw[below] (ra5) node []           {$r_{a_5}$};
\draw[below] (rb3) node []           {$r_{b_3}$};
\draw[below] (rb4) node []           {$r_{b_4}$};
\draw[below] (ra6) node []           {$r_{a_6}$};
\draw[below] (ra7) node []           {$r_{a_7}$};
\draw[below] (rb5) node []           {$r_{b_5}$};

\def\d{4}
\draw (2.5,\d) node [] {$A_2$};
\draw (4.5,\d) node [] {$B_1$};
\draw (7.75,\d) node [] {$A_3$};
\draw (9.97,\d) node [] {$B_2$};
\draw (12.25,\d) node [] {$A_4$};
\draw (15,\d) node [] {$B_3$};
\draw (17.25,\d) node [] {$A_5$};
\draw (19.25,\d) node [] {$B_4$};
\draw (22.75,\d) node [] {$A_6$};
\draw (24.75,\d) node [] {$B_5$};
\draw (27.25,\d) node [] {$A_7$};
\draw (0,\d) node [] {\large $W$};

\draw[left] (a1) node []           {$A_1$};
\draw (1.25,\c+1.5) node []           {$Q_1$};
\draw (6.25,\c+1.5) node []           {$Q_2$};
\draw (11.25,\c+1.5) node []           {$Q_3$};
\draw (16.25,\c+1.5) node []           {$Q_4$};
\draw (21.25,\c+1.5) node []           {$Q_5$};
\draw (26.25,\c+1.5) node []           {$Q_6$};

\draw[dashed, ultra thick] (ra2)--(a2);
\draw[dashed, ultra thick] (ra3)--(a3);
\draw[dashed, ultra thick] (ra4)--(a4);
\draw[dashed, ultra thick] (ra5)--(a5);
\draw[dashed, ultra thick] (ra6)--(a6);
\draw[dashed, ultra thick] (ra7)--(a7);

\draw[dashed, ultra thick] (rb1)--(b1);
\draw[dashed, ultra thick] (rb2)--(b2);
\draw[dashed, ultra thick] (rb3)--(b3);
\draw[dashed, ultra thick] (rb4)--(b4);
\draw[dashed, ultra thick] (rb5)--(b5);

\draw[dashed, ultra thick] (a2) to [bend left = 30] (b2);
\draw[dashed, ultra thick] (a3) to [bend left = 30] (b3);
\draw[dashed, ultra thick] (a4) to [bend left = 30] (b4);
\draw[dashed, ultra thick] (a5) to [bend left = 30] (b5);
\draw[dashed, ultra thick] (a6) to [bend left = 30] (b6);
\draw[dashed, ultra thick] (b1) to [bend right = 15] (a1);
\draw[dashed, ultra thick] (a7) to [bend left = 15] (30,\c+.8);

\draw[dashed, ultra thick] (r1) to (ra2);
\draw[dashed, ultra thick] (rb1) to (ra3);
\draw[dashed, ultra thick] (rb2) to (ra4);
\draw[dashed, ultra thick] (rb3) to (ra5);
\draw[dashed, ultra thick] (rb4) to (ra6);
\draw[dashed, ultra thick] (rb5) to (ra7);
\draw[dashed] (ra2) to (rb1);
\draw[dashed] (ra3) to (ra4);
\draw[dashed] (rb2) to (ra5);
\draw[dashed] (rb3) to (rb4);
\draw[dashed] (ra6) to (ra7);
\draw[dashed] (rb5) to (30,0);

\end{tikzpicture}

\caption{Various paths. The dotted line represents the surface. All the dashed lines represent paths, and the one at the bottom is $R$.
The thick dashed lines are the paths we will use to construct a pair of $S-T$ paths with distance three.
The paths $Q_i$ have been drawn as though their interiors are disjoint from $W$, but they might not be; each is a path of a supercomponent,
but the supercomponent might penetrate into $W$ a short distance, via joints. In particular the ends of each $Q_i$ have been drawn 
exactly
on the surface, but they might be slightly deeper in $W$, or just outside of $W$.}\label{fig:thepaths}
\end{figure}

Now we are ready to construct a pair of $S-T$ paths that are at a distance at least three. (See figure \ref{fig:thepaths}.)
Since $a(H_1) = 0$, there is a vertex $v$ of $H_1$ with $a(v)=0$, and hence with $v\in S$ and $d(v,R)\ge c+4$; let $A_1$ be the one-vertex path consisting of this vertex. For $2\le i\le s$,
let $A_i$ be a shortest
path with
one end in $V(H_i)$ and the other equal to $r_{a_i}$. Similarly, define $B_s$ to be a one-vertex path with vertex some $v\in V(H_s)\cap T$
with $b(v)=0$; 
and for $1\le i<s$ let $B_i$ be a shortest
path with
one end in $V(H_i)$ and the other equal to $r_{b_i}$. It follows that $A_i,B_i$ have length at most $c+1$.
Let $Q_i$ be a path of $H_i$ joining the ends of $A_i$ and $B_i$ in $H_i$.
For $1\le i, j\le n-1$, let $R(i,j)=R(j,i)$ be the subpath of $R$ with ends $r_i,r_j$.
Define $b_0 = 1$ and $a_{s+1}=n-1$.
It follows that the union of the paths 
\begin{itemize}
\item $A_i\cup Q_i\cup B_i$ for $i\in \{1\LL s\}$, odd; and
\item $R(b_i,a_{i+2})$ for $i\in \{0\LL s-1\}$, odd
\end{itemize}
contains an $S-T$ path. (It might not itself be a path, because for instance the paths $B_1,A_3$ might intersect). 
Also, the union of the paths
\begin{itemize}
\item $A_j\cup Q_j\cup B_j$ for $j\in \{1\LL s\}$, even; and
\item $R(b_j,a_{j+2})$ for $j\in \{0\LL s-1\}$, even
\end{itemize}
contains an $S-T$ path.
We must check that these paths have distance at least three. This is true if $s=1$, so we assume that $s\ge 2$.
We have to check that:
\begin{enumerate}
\item[(i)] $d(Q_i,Q_j)\ge 3$ for all $i,j\in \{1\LL s\}$ with $i$ odd and $j$ even;
\item[(ii)] $d(A_i, A_j), d(A_i,B_j), d(B_i,A_j), d(B_i,B_j)\ge 3$ for all $i,j\in \{1\LL s\}$ with $i$ odd and $j$ even;
\item[(iii)] $d(R(b_i,a_{i+2}),R(b_j,a_{j+2}))\ge 3$ for all $i,j\in \{0\LL s-1\}$ with $i$ odd and $j$ even;
\item[(iv)] $d(Q_i, A_j), d(Q_i, B_j), d(A_i, Q_j), d(B_i, Q_j)\ge 3$ for all $i,j\in \{1\LL s\}$ with $i$ odd and $j$ even;
\item[(v)] $d(Q_i, R(b_j,a_{j+2}))\ge 3$ for all $i,j$ with $i\in \{1\LL s\}$, odd and $j\in \{0\LL s-1\}$, even; and 
$d(R(b_i,a_{i+2}), Q_j)\ge 3$ for all $i,j$ with $i\in \{0\LL s-1\}$, odd and $j\in \{1\LL s\}$, even; and
\item[(vi)] $d(A_i,R(b_j,a_{j+2}), d(B_i,R(b_j,a_{j+2})\ge 3$ for all $i,j$ with $i\in \{1\LL s\}$, odd and $j\in \{0\LL s-1\}$, even;
and $d(R(b_j,a_{j+2}), A_j), d(R(b_j,a_{j+2}), B_j)\ge 3$ for all $i,j$ with $i\in \{0\LL s-1\}$, odd and $j\in \{1\LL s\}$, even.
\end{enumerate}

Statement (i) follows from (5). For statement (ii), suppose first that $i,j\ne 1,s$. Since every vertex in $A_i$ has distance at most $c+1$ from $r_{a_i}$,
and the same for $B_i, A_j, B_j$, it suffices to check that $d(r_{a_i},r_{a_j})\ge 2c+5$ and so on. But these four distances are
$|a_i-a_j|,|a_i-b_j|,|b_i-a_j|,|b_i-b_j|$ respectively, and so are all at least $\ell\ge 2c+5$, by (7).
If $i = 1$ and $j\ne s$, the same argument shows that $d(B_1,A_j), d(B_1,B_j)\ge 3$, but we need to check that
$d(A_1, A_j), d(A_1,B_j)\ge 3$. In this case $V(A_1)$ is a vertex $v\in S$, with $d(v,R)\ge c+4$. 
Since $A_j,B_j$ both have length at most $c+1$, it follows that $d(v,A_j),d(v,B_j)\ge 3$. 
The argument is similar if $i=s$ or $j = s$. This proves (ii).

For statement (iii),  note that the subpaths $R(b_i,a_{i+2}),R(b_j,a_{j+2})$ of $R$ are disjoint (by (7)), and the distance 
between them is at least $\ell\ge 3$ (again by (7)). Statement (iv) follows from (6). For statement (v), 
(4) implies that every path from a supercomponent to $R$ has length at least $c-(8-1)/2\ge 3$.
Finally, for statement (vi), we will prove that $d(A_i,R(b_j,a_{j+2})) \ge 3$ for all $i,j$ with $i\in \{1\LL s\}$, odd and 
$j\in \{0\LL s-1\}$, even (the other statement is proved similarly). If $i=1$ then the claim is true since $d(A_1,R)\ge c+4\ge 3$.
So we assume that $i\ge 3$. The distance between $r_{a_i}$ and $R(b_j,a_{j+2})$ is at least $\ell$, by (7); and the distance between 
$r_{a_i}$ and each vertex of $A_i$ is at most $c+1$, and so the distance between $A_i$ and $R(b_j,a_{j+2})$ is at least $\ell-c-1\ge 3$.
This proves \ref{truecase}.~\bbox

\section{Conclusion}\label{conclusion}

While the coarse Menger Conjecture does not hold, there are weakened versions that may be true.  

For $d=2$ (the ``induced paths'' case), the coarse Menger Conjecture remains open for general $k$.  
Gartland, Korhonen and Lokshtanov \cite{gartland} and, independently, Albrechtsen, Huynh, Jacobs, Knappe, and Wollan~\cite{wollan}
have shown that the conjecture holds under an additional bounded degree restriction:

\begin{thm}\label{conjbounded}
{\bf Bounded degree induced Menger:} For every integer $\Delta\ge 1$ there exists $C=C(\Delta)>0$ with the following property.
Let $G$ be a graph, let $k\ge 1$ be an integer, and let $S,T\subseteq V(G)$; then either 
\begin{itemize}
\item there are $k$ paths between $S,T$, pairwise at distance at least two; or 
\item there is a set $X\subseteq V(G)$
with $|X|\le kC$ such that every path between $S,T$ contains a vertex of $X$.
\end{itemize}
\end{thm}

In both cases, the bound on $C$ is quite large (exponential in $\Delta^2$).  However, it is possible that a polynomial bound suffices here, and even that $k-1$ balls of bounded radius suffice to separate $S$ and $T$.  

Bounded degree restrictions seem helpful in finding structures at distance at least two (for example, see also Korhonen \cite{k}).  But for $d>2$, bounding the maximum degree does not help here, as our construction \ref{construction} has maximum degree 3.  
However, a natural weakening of the coarse Menger conjecture is to relax the bound on $|X|$.  Perhaps the following conjecture holds:
\begin{thm}\label{conj3}
{\bf Conjecture:} For all $k,d$ there exist $\ell, n$ such that, if $S,T$ are sets of vertices in a graph $G$, then either
\begin{itemize}
\item there are $k$ paths between $S,T$, pairwise at distance at least $d$; or 
\item there is a set $X\subseteq V(G)$
with $|X|\le n$ such that every path between $S,T$ contains a vertex with distance at most $\ell$ from some member of $X$.
\end{itemize}
\end{thm}

Here we we need $n\ge\lfloor3k/2\rfloor$ (as taking $r$ copies of our construction \ref{construction}, with $S$, $T$ the union of the corresponding sets $S_i$, $T_i$ in each copy, gives a graph which does not contain $2r+1$ paths at pairwise distance three, but no set of $3r-1$ balls of radius at most $\ell$ separate $S$ from $T$). 

Agelos Georgakopoulos brought to our attention another variant of \ref{conj}, which also remains open:
\begin{thm}\label{conj2}
{\bf Conjecture:} For all integers $k,d\ge 1$ there exists $\ell>0$ with the following property.
Let $G$ be a graph and let $S,T\subseteq V(G)$; then either
\begin{itemize}
\item there are $k$ paths between $S,T$, pairwise at distance at least $d$; or
\item there is a vertex $x$ of $G$
such that there do not exist $k-1$ paths between $S$ and $T$, pairwise at distance at least $\ell$ and each with distance at least $\ell$ from $x$.
\end{itemize}
\end{thm}

Finally, we note that our construction has also been used by Davies, Hickingbotham, Illingworth and McCarty \cite{davies} in their disproof of a conjecture of Georgakopoulos and Papasoglu~\cite{agelos} saying (roughly) that graphs without fat $H$-minors are quasi-isometric to graphs without $H$-minors.

\end{document}